\documentclass[12pt]{article}
\usepackage{graphicx}

\usepackage[francais,english]{babel}

\usepackage[]{amsfonts}
\usepackage{amssymb}
\def\couleur(#1 #2 #3)
	{
\special{ps::[end]}}

\def\bx#1{\setbox1=\hbox{\kern3pt{#1}\kern3pt}			
 \dimen1=\ht1 \advance\dimen1 by 3pt \dimen2=\dp1 \advance\dimen2 by 3pt
 \setbox1=\hbox{\vrule height\dimen1 depth\dimen2\box1\vrule}%
 \setbox1=\vbox{\hrule\box1\hrule}%
 \advance\dimen1 by .4pt \ht1=\dimen1
 \advance\dimen2 by .4pt \dp1=\dimen2 \box1\relax}

\def\wbb#1{\kern#1em}
\def\vci{\vrule  width.02em height1.47ex depth-.0ex}		
\def\11{{\rm\wbb{.2}\vci\wbb{-.37}1}}

\newtheorem{Theorem}{Theorem}[section]
\newtheorem{Lemma}[Theorem]{Lemma}
\newtheorem{Definition}[Theorem]{Definition}
\newtheorem{Corollary}[Theorem]{Corollary}

\begin{document}

\title{On separated Carleson sequences in the unit disc of  ${\mathbb{C}}.$ }

\author{Eric Amar}

\date{ }
\maketitle
 \ \par 
\ \par 
\renewcommand{\abstractname}{Abstract}

\begin{abstract}
The interpolating sequences for  $H^{\infty }({\mathbb{D}}),$
  the bounded holomorphic function in the unit disc  ${\mathbb{D}}$
  of the complex plane  ${\mathbb{C}},$  {\small  where characterised
 by L. Carleson by metric conditions on the points. They are
 also characterised by "dual boundedness" conditions which imply
 an infinity of functions. A. Hartmann proved recently that just
 one function in  $H^{\infty }({\mathbb{D}})$  was enough to
 characterize interpolating sequences for  $H^{\infty }({\mathbb{D}}).$
  In this work we  use the "hard" part  of the proof of Carleson
 for the Corona theorem, to extend Hartman's result and answer
 a question he asked in his paper.}\ \par 
\end{abstract}

\section{Introduction.}
Let  ${\mathbb{D}}$  be the unit disc in  ${\mathbb{C}}$  and
  $S$  a sequence of points in  ${\mathbb{D}}.$  Let\ \par 
\quad \quad \quad \quad \quad \quad  $\displaystyle d_{P}(a,b):=\ \left\vert{\frac{a-b}{1-\bar ab}}\right\vert
 $  the pseudo hyperbolic distance and  $d_{H}(a,b):=\tanh ^{-1}(\delta
 (a,b))$  the hyperbolic distance in  ${\mathbb{D}}.$ \ \par 
To say that the sequence  $S$  is separated means that there
 is a  $\eta ,$  such that\ \par 
\quad \quad \quad \quad $\displaystyle \forall a,b\in S,\ a\neq b,\ d_{H}(a,b)\geq \eta
 \iff d_{P}(a,b)\geq \tanh \eta .$ \ \par 
Equivalently to say that the sequence  $S$  is  $\delta $ -separated
 means that the discs  $\forall a\in S,\ D(a,\delta (1-\left\vert{a}\right\vert
 ))$  are disjoint.\ \par 
\quad \quad We shall also need the notion of Carleson measure. Let  $(\zeta
 ,\ h)\in {\mathbb{T}}{\times}(0,1),$  and note\ \par 
\quad \quad \quad \quad \quad \quad  $W(\zeta ,h):=\lbrace z\in {\mathbb{D}}::\left\vert{1-\bar \zeta
 z}\right\vert <h\rbrace $ \ \par 
the associated Carleson window. If  $\nu $  is a borelian measure
 in  ${\mathbb{D}},$  we shall say that  $\nu $  is Carleson
 if there is a constant  $C>0$  such that\ \par 
\quad \quad \quad \quad \quad $\forall \zeta \in {\mathbb{T}},\ \forall h\in (0,1),\ \left\vert{\nu
 }\right\vert (W(\zeta ,h))\leq Ch.$ \ \par 
A sequence  $S$  will be a Carleson sequence if the canonical
 measure associated to it :\ \par 
\quad \quad \quad \quad \quad \quad  $\mu _{S}:=\sum_{a\in S}{(1-\left\vert{a}\right\vert )\delta _{a}},$ \ \par 
is a Carleson measure.\ \par 
\begin{Definition}
We shall say that  $S$  is interpolating  (for  $H^{\infty }({\mathbb{B}})$
 ) if\par 
\quad \quad \quad \quad $\forall \lambda \in l^{\infty }(S),\ \exists f\in H^{\infty
 }({\mathbb{B}})::\forall a\in S,\ f(a)=\lambda _{a}.$ \par 
\end{Definition}
\quad \quad L. Carleson characterized these sequences by the conditions~\cite{CarlInt58}{\bf
  :\ \par 
\quad \quad \quad \quad \quad 	} $\displaystyle (C)\ \ \ \ \ \ \ \ \ \ \ \ \ \ \ \ \ \inf 
 _{a\in S}\prod_{b\in S\backslash \lbrace a\rbrace }{\left\vert{\frac{a-b}{1-\bar
 ba}}\right\vert }>0.$ \ \par 
\quad \quad One can see easily that these conditions are equivalent to the
 fact that  $S$  is {\sl dual bounded }in  $H^{\infty }({\mathbb{D}}),$
  which means :\ \par 
\quad \quad \quad \quad $\exists C>0,\ \forall a\in S,\ \exists \rho _{a}\in H^{\infty
 }({\mathbb{D}}),\ {\left\Vert{\rho _{a}}\right\Vert}_{\infty
 }\leq C::\forall b\in S,\rho _{a}(b)=\delta _{ab}.$ \ \par 
We just take  $\displaystyle \rho _{b}(z):=\frac{B_{b}(z)}{B_{b}(b)}$
  with  $\displaystyle B_{b}(z):=\prod_{a\in S\backslash \lbrace
 b\rbrace }{\frac{a-z}{1-\bar az}\frac{\left\vert{a}\right\vert }{a}}.$ \ \par 
\quad \quad So the {\sl metric conditions (C)} which characterise the interpolation
 are equivalent to the functional characterization namely the
 existence of an infinity of functions verifying the above conditions.\ \par 
\quad \quad \quad 	Another functional characterization is due to D. Sarroste~\cite{D.Amar72}
 :\ \par 
\begin{Theorem}
If there are  $0<\tau <\eta $  such that for any partition  $S=A\uplus
 B$  there is a function  $f\in H^{\infty }({\mathbb{D}}),\ {\left\Vert{f}\right\Vert}_{\infty
 }\leq 1,$  with  $\forall a\in A,\ \left\vert{f(a)}\right\vert
 \leq \tau $  and  $\forall b\in B,\ \left\vert{f(b)}\right\vert
 \geq \eta ,$  then  $S$  is  $H^{\infty }({\mathbb{D}})$  interpolating.\par 
\end{Theorem}
\quad \quad Again one needs an infinity of functions in  $H^{\infty }({\mathbb{D}})$
  to characterize interpolating sequences.\ \par 
\quad \quad A. Hartmann~\cite{Hartmann11} showed that this can be reduced
 to a condition on {\bf only one} function :\ \par 
\begin{Theorem}
Let  $S$  be a separated Blaschke sequence in the unit disc 
 ${\mathbb{D}}$  of  ${\mathbb{C}}.$  There is a partition  $(A,B)$
  of  $S$  such that if there is a function  $f\in H^{\infty
 }({\mathbb{D}})$  with   $f=0$  on  $A$  and  $f=1$  on  $B,$
  then  $S$  is interpolating for  $H^{\infty }({\mathbb{D}}).$ \par 
\end{Theorem}
\quad \quad \quad 	So a natural question after these results is : is it possible
 to have an analogous result as D. Sarroste's one but replacing
 {\sl for any partition} by {\sl there is a partition} ?\ \par 
\quad \quad 	The aim of this work is to prove that the answer is yes, provided
 that  $\tau <\eta ^{\kappa }$  for a certain constant  $\kappa
 >1$  introduced by Carleson in his proof of the corona theorem.\ \par 
\quad \quad \quad 	We shall need the following notions.\ \par 
\begin{Definition}
We shall say that the partition  $(A,B)$  of the sequence of
 points  $S\subset {\mathbb{D}}$  is "good" if there is  $\varphi
 \ :\ A\rightarrow B$  such that   $\forall a\in A,\ d_{H}(a,\varphi
 (a))=\inf  _{c\in S\backslash \lbrace a\rbrace }d_{H}(a,c)$
  and if there is  $\psi \ :\ B\rightarrow A$  such that   $\forall
 b\in B,\ d_{H}(b,\psi (b))=\inf  _{c\in S\backslash \lbrace
 b\rbrace }d_{H}(b,c).$ \par 
\end{Definition}
\quad \quad We shall need more specific partitions. 	Let  $\gamma \in \rbrack
 0,1\lbrack $  ; we set\ \par 
\quad \quad \quad \quad \quad \quad  $\displaystyle C_{n}=C_{n}(\gamma ):=\lbrace z\in {\mathbb{D}}::1-\gamma
 ^{n+1}<\left\vert{z}\right\vert \leq 1-\gamma ^{n}\rbrace .$ \ \par 

\subsubsection{Restricted good partition.}
\begin{Definition}
A restricted good partition of the discrete sequence  $S$  in
 the disc is a partition  $(A,B)$  of  $S$  such that 	 $\exists
 \gamma \in \rbrack 0,1\lbrack $  and with  $A_{n}:=A\cap C_{n}(\gamma
 ),\ B_{n}:=B\cap C_{n}(\gamma )$  we have  $(A_{n},B_{n})$ 
 is a good partition of  $S_{n}:=S\cap C_{n}(\gamma )$  for any
  $n\in {\mathbb{N}}.$ \par 
\end{Definition}
\quad \quad \quad 	As we shall see later there are always restricted good partition
 for a discrete sequence  $S$  in the disc.\ \par 

\subsubsection{Hoffman partition.}
Let  $S$  be a discrete sequence in the disc. (See J. Garnett~\cite{garnett}.)\
 \par 
\quad \quad \quad 	We shall cut  $S$  in two parts ; for this let\ \par 
\quad \quad \quad \quad \quad \quad  $D_{1}:=\lbrace z\in {\mathbb{D}}::\mathrm{A}\mathrm{r}\mathrm{g}z\in
 \lbrack 0,\pi \lbrack \rbrace ,\ D_{2}:=\lbrace z\in {\mathbb{D}}::\mathrm{A}\mathrm{r}\mathrm{g}z\in
 \lbrack \pi ,\ 2\pi \lbrack \rbrace .$ \ \par 
Now set  $S_{1}:=S\cap D_{1},\ S_{2}:=S\cap D_{2}.$  Because
 if  $S_{1}$  and  $S_{2}$  are Carleson sequences then  $S$
  is also a Carleson sequence, it will be enough to deal with
 one of them, say  $S=S_{1}.$ \ \par 
\quad \quad 	We start with the point  $a_{0}=\left\vert{a_{0}}\right\vert
 e^{i\theta _{0}}$  in  $S_{n}:=C_{n}\cap S$  with the smallest
 module and smallest argument ; if  $\# S_{n}=1,$  put  $a_{0}$
  in  $A_{n}$  and set  $B_{n}:=\emptyset \ ;$  if not take the
 next point in  $S_{n}$  with the same argument as  $a,$  hence
 with a bigger modulus, if any, or at the right side of  $a_{0},$
  i.e. such that its argument  $\theta $  is bigger than  $\theta
 _{0}.$  Call it  $b_{0}$  and define  $\varphi (a_{0}):=b_{0}.$ \ \par 
Now take the next point in  $S_{n}$  at the right side of  $b_{0},$
  i.e. the same way as above, and call it  $a_{1}$  etc... Then
 each time define  $\displaystyle \varphi (a_{j}):=b_{j}.$  Call
  $A_{n}$  the set of all  $a_{j}'s$  and  $B_{n}$  the set of
 all  $b_{j}'s.$   Because  $S$  is discrete  $A_{n}$  and  $B_{n}$
  are finite.\ \par 
\begin{Definition}
Let  $A:=\bigcup_{n\in {\mathbb{N}}}{A_{n}}$  and  $B:=\bigcup_{n\in
 {\mathbb{N}}}{B_{n}}$  then  $(A,B)$  is the Hoffman partition of  $S.$ \par 
\end{Definition}
\quad \quad \quad 	We see easily that for any  $a\in A,\ \varphi (a)$  has always
 its argument bigger or equal to the argument of  $a.$ \ \par 
\begin{Definition}
Let  $(A,B)$  be a restricted good partition or a Hoffman partition
 of the sequence  $S\subset {\mathbb{D}}.$  Let  $\kappa \geq
 1$  be a constant, the sequence  $S\subset {\mathbb{D}}$  is
  $\kappa $ -ultra-separated if  $S$  is separated and if there
 are  $0<\tau <\eta ,\ \tau <\eta ^{\kappa }$  and a function
  $f$  in  $H^{\infty }({\mathbb{D}}),\ {\left\Vert{f}\right\Vert}_{\infty
 }\leq 1,$  such that  $\ \left\vert{f}\right\vert \leq \tau
 $  on  $A$  and  $\ \left\vert{f}\right\vert \geq \eta $  on  $B.$ \par 
\end{Definition}
\quad \quad Now we can state the theorem.\ \par 
\begin{Theorem}
\quad ~\label{ultraSepBall61}There is a constant  $\kappa >1$  such
 that the sequence  $S$  is  $H^{\infty }({\mathbb{D}})$  interpolating
 if and only if it is  $\kappa $ -ultra-separated.\par 
\end{Theorem}
\quad \quad \quad 	This constant  $\kappa $  was introduced by Carleson in his
 proof of the corona theorem.\ \par 
The theorem~\ref{ultraSepBall61} generalizes the result of A.
 Hartman and answer {\sl positively} to his question :\ \par 
\quad \quad \quad 	if there is a  $f\in H^{\infty }({\mathbb{D}})$  such that 
 $\forall a\in A,\ f(a)=0,\ \forall b\in B,\ \left\vert{f(b)}\right\vert
 \geq \eta >0$  for a Hoffman partition  $(A,B)$  of  $S,$  and
  $S$  separated, is  $S$  interpolating ?\ \par 
\ \par 
\quad \quad I introduce good partitions for dealing with this problem in
 the unit ball of  ${\mathbb{C}}^{n}.$  This notion in invariant
 by automorphisms and hence natural. The result in the ball,
 not as good as in the disc, will be posted later. It involves
 complex geometry and the key fact is that the measure  $(1-\left\vert{z}\right\vert
 )\left\vert{\partial f(z)}\right\vert ^{2}dm(z)$  is a Carleson
 measure in the unit ball of  ${\mathbb{C}}^{n}.$ \ \par 

\section{General facts.}
\begin{Lemma}
~\label{ultraSepBall135}Let  $S$  be a discrete sequence in the
 metric space  $(X,\ d),$  there is a good partition  $(S_{1},\
 S_{2})$  of  $S.$ \par 
\end{Lemma}
\quad \quad Proof.\ \par 
Take a point  $O\in X$  and  $a_{1}\in S$  such that  $d(a_{1},\
 O)$  is minimal, if  $\# S=1$  set  $\varphi (a_{1})=a_{1}$
  and  $S_{1}:=\lbrace a_{1}\rbrace =S\ ;\ S_{2}=\emptyset \
 ;$  then we are done.\ \par 
If  $\# S\geq 2,$  then take  $b_{1}\in S$  a nearest neighbour
 for the distance  $d$  of  $a_{1}$  and define  $\varphi (a_{1})=b_{1}.$
  By the assumption on the cardinality of  $S,\ b_{1}$  exists.
 Take  $a_{2}$  a nearest neighbour of  $b_{1},$  if it exists,
 and define   $\psi (b_{1}):=a_{2}\ ;$  if  $a_{2}=a_{1}$  we
 stop at this "perfect" pair  $(a_{1},b_{1})$  with  $\psi (b_{1}):=a_{1}.$
  If not we continue with  $b_{2}$  nearest neighbour of  $a_{2}$
  etc... We stop at a perfect pair. This way we get a branch
   $B_{1}$  finite or infinite. We put all the  $"a"$  in  $S_{1}$
  and all the  $"b"$  in  $S_{2}.$ \ \par 
If it remains points in  $S$  we have that the points in   $S\backslash
 B_{1}$  are far from the points in  $B_{1}$  by construction.
 We take a point  $c$  in  $S\backslash B_{1}$  the nearest from  $O.$ \ \par 
\quad \quad A) If all the nearest points from  $c$  are in  $B_{1},$  which
 may happen, we take one of them,  $d,$  now if  $d$  is in 
 $S_{1},$  we put  $c$  in  $S_{2}$  and we set  $\psi (c):=d.$
  If  $d$  is in  $S_{2},$  we put  $c$  in  $S_{1}$  and we
 set  $\varphi (c):=d.$  This completes  $B_{1}$  and we start all again.\ \par 
\quad \quad B) If  $c$  has a nearest neighbour which is not in  $B_{1},$
  we start a new branch  $B_{2}$  etc... \ \par 
A new point may have its nearest neighbour in  $B_{1}$  or in
  $B_{2},$  etc... Then  we put it in  $B_{1}$  or in  $B_{2},\
 ...$  as in the step A.\ \par 
\quad \quad We continue this way in order to exhaust  $S.$ \ \par 
\quad \quad The   $S_{1}$  part is all the  $"a"$  and  $S_{2}$  is all the  $"b".$ \ \par 
Then  $S$  is a bipartite graph with components  $S_{j},\ j=1,2$
  on which the two applications  $\varphi ,\ \psi $  are well
 defined.  $\hfill\blacksquare $ \ \par 

\section{Proof of the main theorem.}
\quad \quad \quad 	Let  $S$  be a discrete sequence in the unit disc  ${\mathbb{D}}.$
  Fix any  $0<\gamma <1$  and  $n\in {\mathbb{N}}$  and recall\ \par 
\quad \quad \quad \quad \quad \quad  $C_{n}=C_{n}(\gamma ):=\lbrace z\in {\mathbb{D}}::1-\gamma ^{n+1}<\left\vert{z}\right\vert
 \leq 1-\gamma ^{n}\rbrace ,\ S_{n}:=S\cap C_{n}.$ \ \par 
We shall use the following lemmas.\ \par 
\begin{Lemma}
~\label{SepDisc59}If the number of points in  $S_{n}$  is smaller
 than  $m$ , i.e.  $\ \forall a\in S,\ \# S_{n}\leq m,$  then
  $S$  is a Carleson sequence.\par 
\end{Lemma}
\quad \quad \quad 	Proof.\ \par 
Let  $W=W(\zeta ,h)$  be a Carleson window ; for  $S$  to be
 Carleson we must have\ \par 
\quad \quad \quad \quad \quad \quad  $\displaystyle \exists C::\sum_{a\in S\cap W}{(1-\left\vert{a}\right\vert
 )}\leq Ch.$ \ \par 
If  $a\in W$  we have  $1-\left\vert{a}\right\vert \leq h,$  hence\ \par 
\quad \quad \quad \quad \quad \quad  $\displaystyle \ \sum_{a\in S\cap W}{(1-\left\vert{a}\right\vert
 )}\leq \sum_{a\in S,\ 1-\left\vert{a}\right\vert \leq h}{(1-\left\vert{a}\right\vert
 )}.$ \ \par 
But  $a\in C_{n}(\gamma )\Rightarrow \gamma ^{n+1}<1-\left\vert{a}\right\vert
 \leq \gamma ^{n},$  and because there are at most  $m$  points
 in  $S\cap C_{n}(\gamma )$  we have\ \par 
\quad \quad \quad \quad \quad \quad  $\displaystyle \ \sum_{a\in S,\ 1-\left\vert{a}\right\vert \leq
 h}{(1-\left\vert{a}\right\vert )}=\sum_{n::\gamma ^{n+1}<h}{\sum_{a\in
 C_{n}(\gamma )\cap S}{(1-\left\vert{a}\right\vert )}}\leq m\sum_{n::\gamma
 ^{n+1}<h}{\gamma ^{n}}\leq \frac{m}{\gamma (1-\gamma )}h.$ \ \par 
\quad \quad \quad 	Hence we have the lemma with  $C=m/\gamma (1-\gamma ).$   $\hfill\blacksquare
 $ \ \par 
\begin{Lemma}
Let  $S$  be a discrete sequence in  ${\mathbb{D}},$  then there
 is a restricted good partition for  $S.$ \par 
\end{Lemma}
\quad \quad \quad 	Proof.\ \par 
Take any  $\gamma \in \rbrack 0,1\lbrack \ ;$  because  $S$ 
 is discrete,  $S_{n}:=S\cap C_{n}(\gamma )$  has only a finite
 number of points. If  $S_{n}=\emptyset ,$  we simply set  $A_{n}=B_{n}:=\emptyset
 .$  If its cardinal is bigger than one, we can apply the general
 lemma~\ref{ultraSepBall135} : there is a good partition  $(A_{n},B_{n})$
  of  $S_{n}.$ \ \par 
Setting  $A=\bigcup_{n\in {\mathbb{N}}}{A_{n}},\ B=\bigcup_{n\in
 {\mathbb{N}}}{B_{n}},$  we have that  $(A,B)$  is a restricted
 good partition of  $S.$   $\hfill\blacksquare $ \ \par 
We remark that  $\displaystyle \gamma \leq \frac{1-\left\vert{a}\right\vert
 }{1-\left\vert{\varphi (a)}\right\vert }\leq 1/\gamma $  because
  $a$  and  $\varphi (a)$  belong to the same  $C_{n}(\gamma ).$ \ \par 
\ \par 
\quad \quad \quad 	Back to the proof of the main theorem.\ \par 
Let  $W=W(\zeta ,h)$  be a Carleson window ; we have to show
 that  $S$  is Carleson i.e.\ \par 
\quad \quad \quad  $\displaystyle \ \sum_{a\in A\cap W}{(1-\left\vert{a}\right\vert
 )}\leq Ch,$ \ \par 
then, because  $S$  is separated, it will be  $H^{\infty }({\mathbb{D}})$
  interpolating.\ \par 
\quad \quad \quad 	We shall cut the set  $A\cap W$  in two parts.\ \par 

\section{The points  $E_{W}:=a\in W\cap S$  such that  $\varphi (a)\notin W.$ }

\subsection{Case of Hoffman partition.}
\quad \quad \quad 	We shall work in the half plane  ${\mathbb{C}}^{+},$  because
 the geometry is easiest. This means that  $1-\left\vert{a}\right\vert
 $  is replaced by  $\Im a,$  the imaginary part of  $a,$  and
  $\mathrm{A}\mathrm{r}\mathrm{g}a$  is replaced by  $\Re a,$
  the real part of  $a.$  Let  $C_{a}$  be the strip  $C_{n}(\gamma
 )=\lbrace z::\gamma ^{n+1}<\Im z\leq \gamma ^{n}\rbrace $  which
 contains  $a.$  If we deal with a Hoffman partition, the point
  $b:=\varphi (a)$  is the nearest point in  $S\cap C_{a},$ 
 either with the same real part, hence with a smaller imaginary
 part, or on the right of  $a$  ; hence if  $b\notin W$  this
 means that there is no points of  $A$  between  $a$  and the
 right vertical side of  $W$  in the strip  $C_{n}(\gamma )$
  containing  $a.$  So  $a$  is the nearest point in  $A$  to
 the right side of  $W$  in the strip  $C_{a}.$  We shall take
 the maximum possible of these points which means that we have
 at most one point in each  $C_{n}(\gamma )\cap W$  and then we have\ \par 
\quad \quad \quad \quad \quad \quad  $\ \sum_{a\in E_{W}}{\Im a}\leq Ch,$ \ \par 
by lemma~\ref{SepDisc59}.\ \par 

\subsection{Case of restricted good partition.}
\quad \quad We shall work again in the half plane because the geometry is
 easiest. So  $W$  is a square with one side on the real axis.
 Let  $c$  be the orthogonal projection of  $a$  on the side
 of  $W$  in the direction of  $b.$ \ \par 
We define the {\sl border strip} to be a tube  $T(a)$  around
 the segment  $\lbrack a,c\rbrack $  of width  $r\Im a.$ \ \par 
\ \par 
\quad \quad \quad 	The partition  $(A,B)$  being restricted, this means that  $b$
  belongs to the same strip   $C_{n}(\gamma ):=\lbrace z\in {\mathbb{C}}^{+}::\gamma
 ^{n+1}<\Im z\leq \gamma ^{n}\rbrace $  as  $a.$  Let us denote
  $C_{a}$  the strip  $C_{n}(\gamma )$  to which  $a$  belongs.\ \par 
\begin{Lemma}
~\label{SepDisc48}Let  $(A,B)$  be a restricted good partition
 of  $S$  in  ${\mathbb{D}}.$   Let  $W=W(\zeta ,h)$  be a Carleson
 window and  $a\in A$  and  $b:=\varphi (a)$  be such that  $a\in
 W,\ b\notin W.$  Then the border strip  $T(a)$  contains at
 most a fixed number  $m$  of points of  $A.$ \par 
\end{Lemma}
\quad \quad \quad 	Proof.\ \par 
Because  $b$  is the nearest point in  $S\cap C_{a}$  we have
 that there is no point of  $S\cap C_{a}$  in the hyperbolic
 ball  $Q(a,b)$  "centered" at  $a$  and passing through  $b.$
  So the worst case is when  $b$  belongs to one of the three
 sides of  $W$  in  ${\mathbb{C}}^{+}.$  Suppose first that 
 $b$  is in the vertical left side of  $W.$ \ \par 
We have the worst case when  $C_{a}=C_{n}(\gamma ),$  with  $\Im
 a=\gamma ^{n+1}.$  Then the border strip is\ \par 
\quad \quad \quad \quad \quad \quad  $T(a):=\lbrace z=x+iy\in W::(1-r)\Im a<y<(1+r)\Im a,\ x<\Re a\rbrace .$ \ \par 
See the picture below.\ \par 
\includegraphics[trim=0cm 6cm 0cm 6cm, clip, width=6cm]{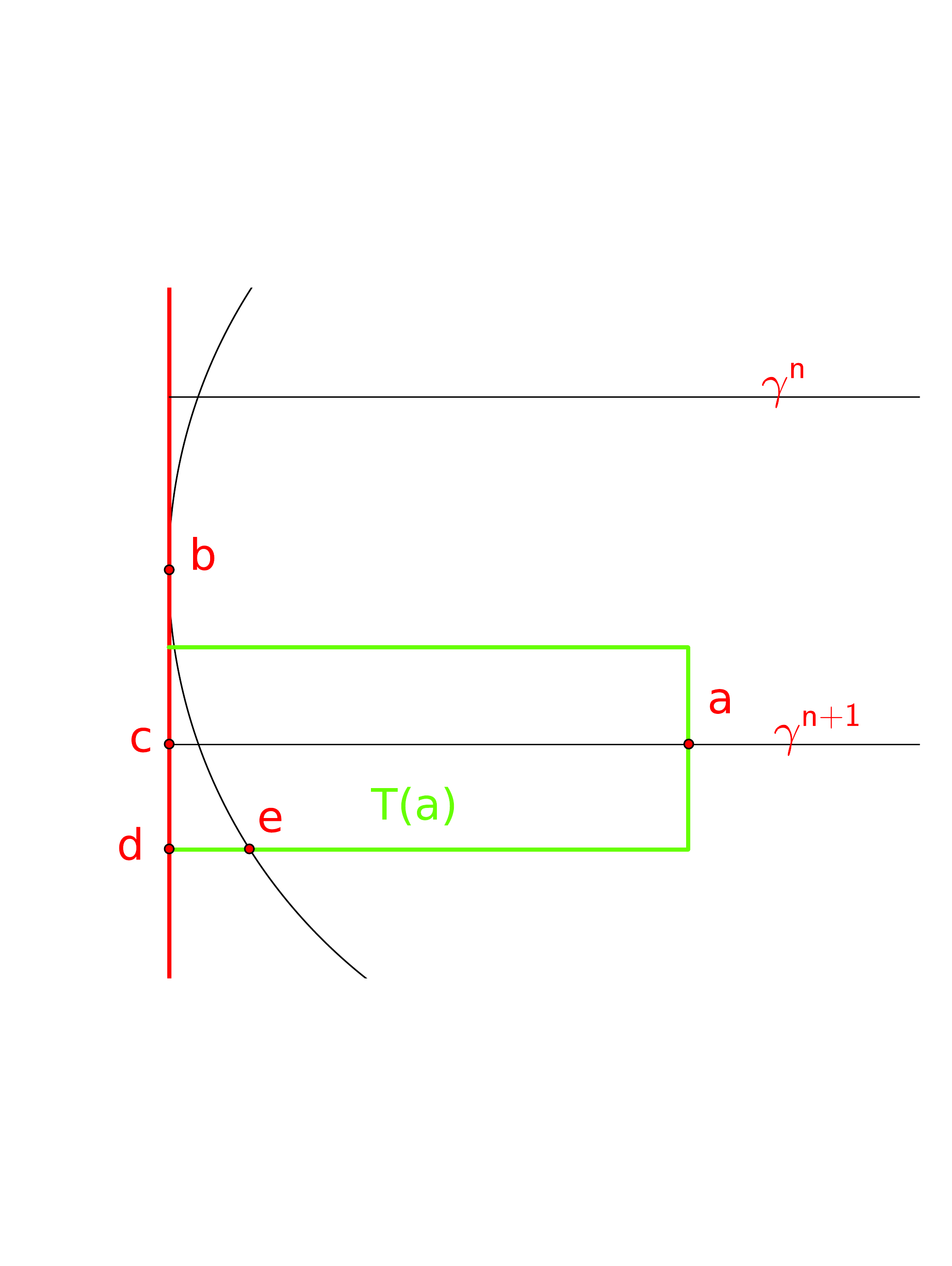}
\
 \par 
\quad \quad \quad 	Then points of  $A\backslash \lbrace a\rbrace $  in  $T(a)$
  must be in the triangle  $bde\cap T(a)$  but  $\ \left\vert{de}\right\vert
 \leq \left\vert{bd}\right\vert \leq \gamma ^{n}(1-\gamma )$
  and in this triangle there is at most  $m=m(\gamma ,\delta
 )$  points in  $A$  because  $A$  is a  $\delta $ -separated
 sequence. A fortiori the number of points in  $A$  such that
  $b\notin W$  is smaller than  $m.$ \ \par 
\quad \quad \quad 	If  $b$  is on the right side of  $W,$  this is the same.\ \par 
\ \par 
\quad \quad \quad 	Suppose now that  $b$  is on the top of  $W.$ \ \par 
\includegraphics[trim=0cm 7cm 0cm 7cm, clip, width=7cm]{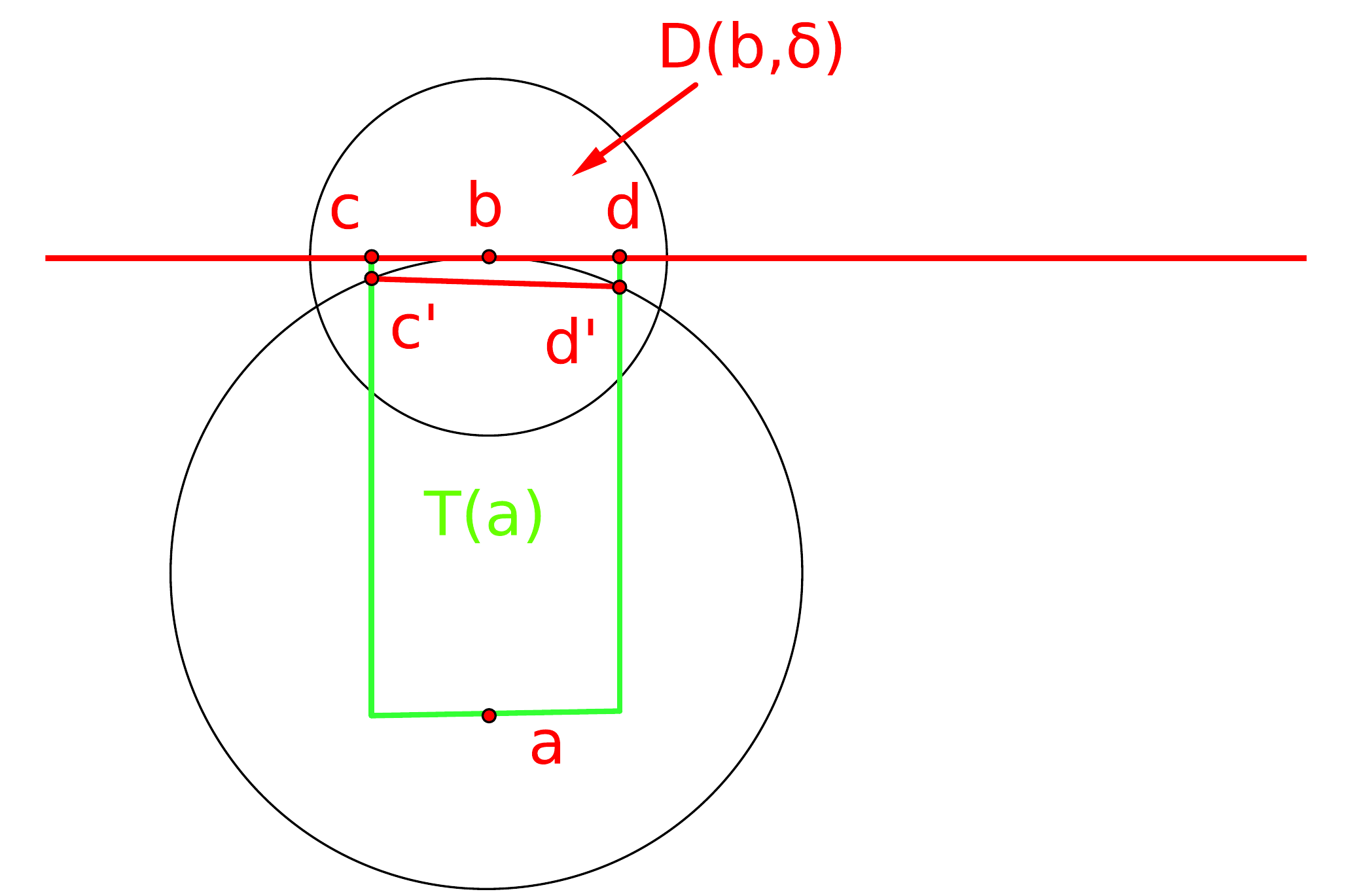}
\
 \par 
Still because  $b$  is the nearest point in  $S\cap C_{a}$  to
  $a,$  the points in  $(A\backslash \lbrace a\rbrace )\cap T(a)\cap
 W$  must be in a rectangle  $cdc'd'$  of sides less than  $2r\Im
 a.$  Because  $S$  is  $\delta $ -separated, there is no point
 of  $S\backslash \lbrace b\rbrace $  in the disc  $D(b,\delta
 (1-\left\vert{b}\right\vert )),$  hence there is no point of
  $S$  in this rectangle provided that  $r<\delta ,$  but the
 point  $b,$  and the lemma.  $\hfill\blacksquare $ \ \par 
\quad \quad Again taking border strips  $T(a)$  with half of the width, they
 become disjoint. Now set  $E_{W}'$  the points in  $E_{W}$ 
 such that the  $T(a)$  is based on the vertical sides of  $W$
  and  $E_{W}''$  the points in  $E_{W}$  such that the  $T(a)$
  is based on the top of  $W.$ \ \par 
We have, by lemma~\ref{SepDisc59}, that  $\displaystyle \ \sum_{a\in
 E_{W}'}{\Im a}\lesssim h.$ \ \par 
For  $E_{W}''$  we have that  $a\in E_{W}''$  must belong to
 the  $C_{n}(\gamma )$  with the smallest  $n$  such that  $W\cap
 C_{n}(\gamma )\neq \emptyset .$ So we have  $\gamma ^{n}\leq
 \Im a\leq \gamma ^{n-1},$  hence, because the tube  $T(a)$ 
 has width  $r\Im a,$  we have at most  $h/r\gamma ^{n}$  such
 tubes. So finally\ \par 
\quad \quad \quad \quad \quad \quad  $\displaystyle \ \sum_{a\in E_{W}''}{\Im a}\leq \gamma ^{n-1}{\times}\frac{h}{r\gamma
 ^{n}}\leq \frac{1}{r\gamma }h.$ \ \par 
So adding these two inequalities we get\ \par 
\quad \quad \quad \quad \quad \quad  $\displaystyle \ \sum_{a\in E_{W}}{\Im a}\lesssim h,$ \ \par 
and the right estimate for  $E_{W}.$ \ \par 

\section{The points  $F_{W}:=a\in W\cap S$  such that  $\varphi (a)\in W.$ }
\quad \quad In this part, we shall use the "hard" part of the proof of L.
 Carleson of the corona theorem, as  interpreted by H\"ormander~\cite{HormCor67}
 (Lemma 11, p 948):\ \par 
\begin{Lemma}
~\label{ultraSepBall62}There exists a constant  $\kappa $  such
 that if  $\displaystyle 0<\eta <\frac{1}{2}$  and  $f\in H^{\infty
 }({\mathbb{D}}),\ \sup \left\vert{f}\right\vert \leq 1,$  one
 can find  $\psi $  with  $0\leq \psi \leq 1$  so that  $\displaystyle
 \ \frac{\partial \psi }{\partial \bar z}dm$  is a Carleson measure
 in  ${\mathbb{D}}$  and\par 
\quad \quad \quad \quad \quad \quad  $\psi (z)=0$  when   $\ \left\vert{f(z)}\right\vert <\eta ^{\kappa
 },\ \psi (z)=1$  when  $\ \left\vert{f(z)}\right\vert \geq \eta .$ \par 
\end{Lemma}
We shall call this  $\kappa $  the Carleson constant. Because
  $\psi $  is real valued we also have  $\displaystyle \ \frac{\partial
 \psi }{\partial z}dm$  is Carleson hence  $\ \left\vert{{\rm{grad
 }}\psi }\right\vert dm$  is Carleson.\ \par 
\ \par 
\quad \quad Let  $(A,B)$  be the partition of  $S$  associated to the function
  $f\in H^{\infty }({\mathbb{D}}).$  Taking eventually a power
 of  $f,$  we can assume that  $\eta <1/2$  to fit with the hypotheses
 of H\"ormander's lemma.\ \par 
\ \par 
We shall use the following well known facts :\ \par 
\quad \quad \quad 	1- If  $f\in H^{\infty }({\mathbb{D}})$  and  $\ \left\vert{f(a)}\right\vert
 \leq \tau $  then if  $\tau '>\tau $  there is a  $r>0$  depending
 only on  $f,\tau '$  such that :\ \par 
\quad \quad \quad \quad \quad \quad  $\forall z\in D(a,r(1-\left\vert{a}\right\vert )),\ \left\vert{f(z)}\right\vert
 <\tau '.$ \ \par 
\quad \quad \quad 	2- If  $f\in H^{\infty }({\mathbb{D}})$  and  $\ \left\vert{f(b)}\right\vert
 \geq \eta $  then if  $\eta '<\eta $  there is a  $r>0$  depending
 only on  $f,\eta '$  such that :\ \par 
\quad \quad \quad \quad \quad \quad  $\forall z\in D(b,r(1-\left\vert{b}\right\vert )),\ \left\vert{f(z)}\right\vert
 >\eta '.$ \ \par 
\ \par 
Let again  $(A,B)$  the restricted or Hoffman partition of  $S$
  and  $f\in H^{\infty }({\mathbb{D}})$  the function with\ \par 
\quad \quad \quad \quad \quad \quad  $\forall a\in A,\ \left\vert{f(a)}\right\vert \leq \tau <\eta
 ^{\kappa },\ \forall b\in B,\ \left\vert{f(b)}\right\vert \geq \eta .$ \ \par 
By the fact 1 we have\ \par 
\quad \quad \quad \quad \quad \quad  $\displaystyle \forall z\in D(a,r_{1}(1-\left\vert{a}\right\vert
 )),\ \left\vert{f(z)}\right\vert <\tau '$ \ \par 
and by fact 2\ \par 
\quad \quad \quad \quad \quad \quad  $\displaystyle \forall z\in D(b,r_{2}(1-\left\vert{b}\right\vert
 )),\ \left\vert{f(z)}\right\vert >\eta '.$ \ \par 
Take  $r=\min (r_{1},r_{2})$  to have both.\ \par 
\ \par 
\quad \quad We shall need the following notions.  	Let  $0<\alpha <1,\ \beta
 >0,\ a,b\in {\mathbb{D}}$  and set  $R(a,b,\alpha ,\beta )$
  a tube around a smooth curve  $\Gamma (a,b)$  with thickness
  $\tau :=\alpha \min (1-\left\vert{a}\right\vert ,1-\left\vert{b}\right\vert
 )$  i.e.\ \par 
\quad \quad \quad \quad \quad \quad  $R(a,b,\alpha ,\beta ):=\bigcup_{c\in \Gamma (a,b)}{D(c,\tau )}.$ \ \par 
Moreover the Lebesgue measure on  $\displaystyle R(a,b,\alpha
 ,\beta )$  must be  $\beta $ -equivalent to the Lebesgue measure
 on the product  $(-\tau ,\tau ){\times}\Gamma (a,b).$ \ \par 
\begin{Lemma}
~\label{SepDisc610}Let  $(A,B)$  be a restricted or a Hoffman
 partition of  $S.$  Then we can make tubes  $\displaystyle R(a,\varphi
 (a),s,\pi /2)$  which are disjoint.\par 
\end{Lemma}
\quad \quad \quad 	Proof.\ \par 
If  $(A,B)$  is a restricted good partition, then for any  $a\in
 A$  we take the tube of width  $r(1-\left\vert{a}\right\vert
 )$  around the segment  $(a,\varphi (a)).$  Because  $a,\ b:=\varphi
 (a)$  belong to the same strip  $C_{n}(\gamma ),$  we have that
  $\displaystyle \gamma \leq \frac{1-\left\vert{a}\right\vert
 }{1-\left\vert{\varphi (a)}\right\vert }\leq 1/\gamma .$  Also
  $S$  being  $\delta $ -separated, there is no point of  $S$
  in a disc  $D(b,\delta (1-\left\vert{b}\right\vert ))$  centered
 at  $b$  and of radius  $\delta (1-\left\vert{b}\right\vert
 ).$  Moreover, because  $b$  is the nearest point in  $S$  to
  $a,$  in the hyperbolic metric, the disc  $Q(a,b)$  of all
 the points in  ${\mathbb{D}}$  nearer to  $a$  than  $b$  contains
 no other points of  $S$  ; hence the tube with width  $\displaystyle
 s=\min (\delta ,\ r)$  does not contains any point in  $S$ 
 but  $a$  and  $b.$  See the picture :\ \par 
\includegraphics[trim=0cm 7cm 0cm 7cm, clip, width=7cm]{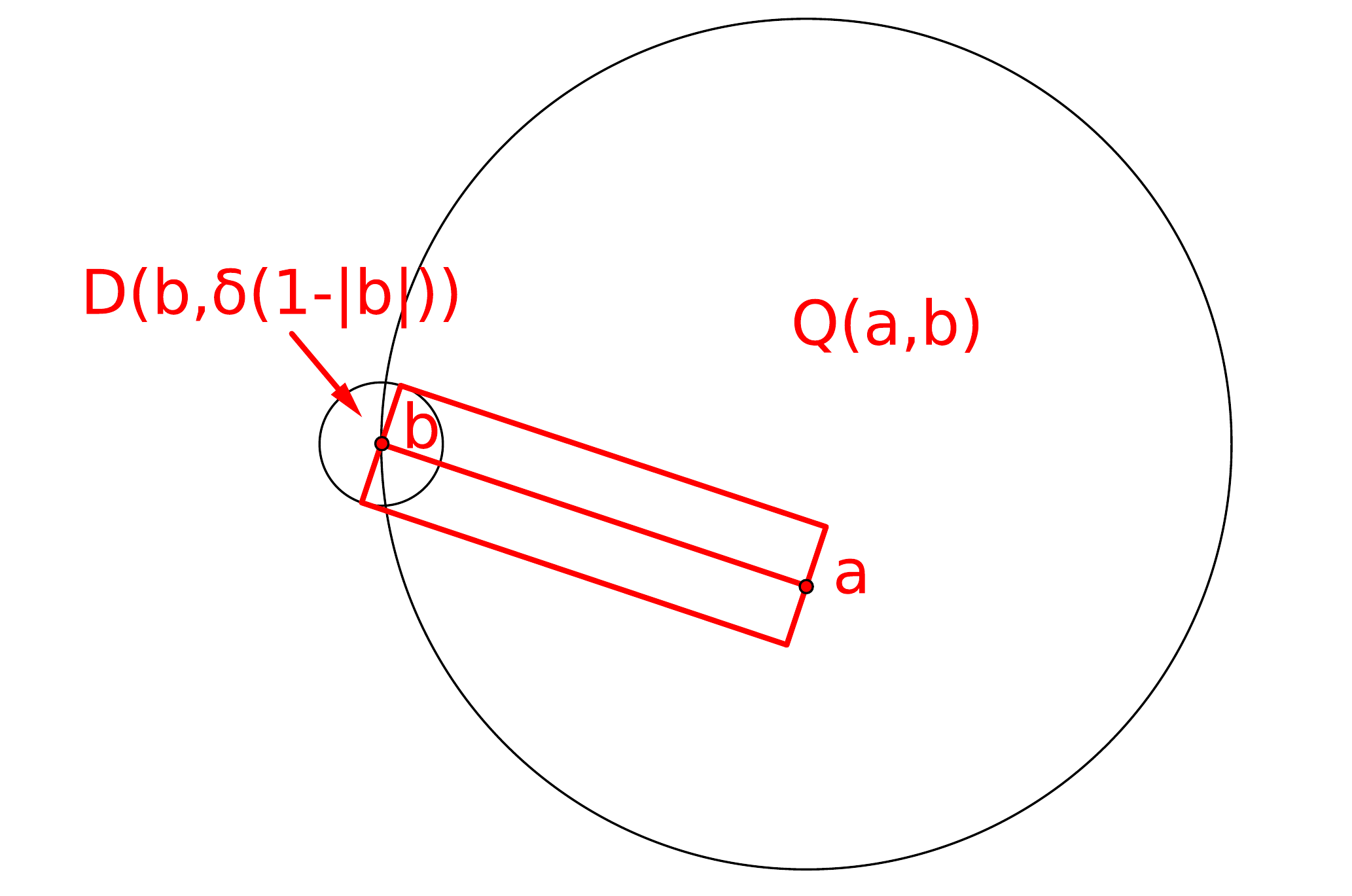}
\
 \par 
\quad \quad \quad 	So we take the tubes with width  $s/2$  instead of  $s$  to
 have disjoint tubes.\ \par 
\ \par 
If  $(A,B)$  is a Hoffman partition, again  $a$  and  $b:=\varphi
 (a)$  are in the same strip  $C_{n}(\gamma )$  hence  $\displaystyle
 1-\left\vert{a}\right\vert \simeq 1-\left\vert{b}\right\vert .$ \ \par 
Because again the points in  $S$  are  $\delta $ -separated,
 we can perturb a little bit the segment  $(a,b)$  in order to
 avoid discs  $D(c,\delta (1-\left\vert{c}\right\vert )),\ c\in
 S\backslash \lbrace a,b\rbrace $  with a curve whose length
 is less than  $\pi /2$  times the length of  $(a,b).$ \ \par 
\includegraphics[trim=0cm 7cm 0cm 7cm, clip, width=7cm]{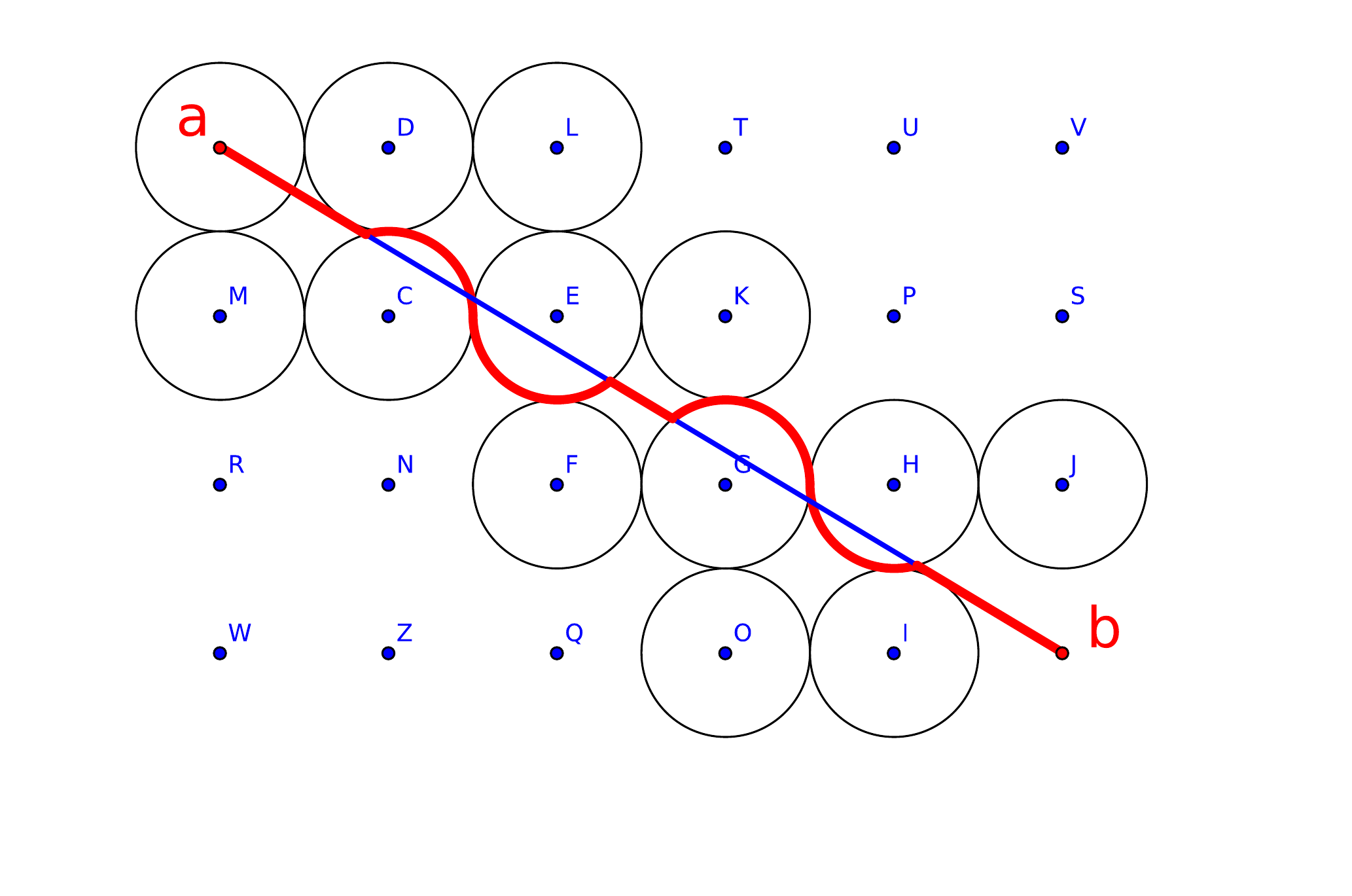}
\
 \par 
\quad \quad \quad 	This means that around this curve we can make a tube of width
 less than  $(\delta /4){\times}(1-\left\vert{a}\right\vert )$
  and these tubes are still disjoint provided that no segments
  $(a,\varphi (a))$  and  $(a',\varphi (a'))$  cut each other.
 But by the construction of the Hoffman partition this cannot
 happen. We have now to take	  $s=\min (\delta /4,\ r)$  to have
 that the tubes  $\displaystyle R(a,\varphi (a),s,\ \pi /2)$
  are disjoint.\ \par 
Hence the lemma.  $\hfill\blacksquare $ \ \par 
\ \par 
Because  $S$  is  $\kappa $ -ultra-separated, we have\ \par 
\quad \quad \quad \quad \quad \quad  $\exists f\in H^{\infty }({\mathbb{D}})::{\left\Vert{f}\right\Vert}_{\infty
 }\leq 1,\ \forall a\in A,\ \left\vert{f(a)}\right\vert \leq
 \tau <\eta ^{\kappa },\ \forall b\in B,\ \left\vert{f(b)}\right\vert
 \geq \eta .$ \ \par 
Fix  $\displaystyle \tau '::\tau <\tau '<\eta ^{\kappa }$  then
 we have a  $r>0$  such that\ \par 
\quad \quad \quad \quad \quad \quad  $\forall a\in A,\ \forall z\in D(a,r(1-\left\vert{a}\right\vert
 )),\ \left\vert{f(z)}\right\vert \leq \tau '<\eta ^{\kappa }.$  \ \par 
By lemma~\ref{SepDisc610} we have that the length of the curve
  $\Gamma (a,\varphi (a))$  is smaller than  $\pi /2$  times
 the length of the segment  $(a,\varphi (a))$  so we can enlarge
 a little bit  $W$  say  $W':=W'(\zeta ,\pi h/2),$  in order to have that\ \par 
\quad \quad \quad \quad \quad \quad  $\displaystyle R(a,\varphi (a),s,\pi /2)\subset W'.$ \ \par 
See the picture\ \par 
\includegraphics[trim=0cm 8cm 0cm 8cm, clip, width=8cm]{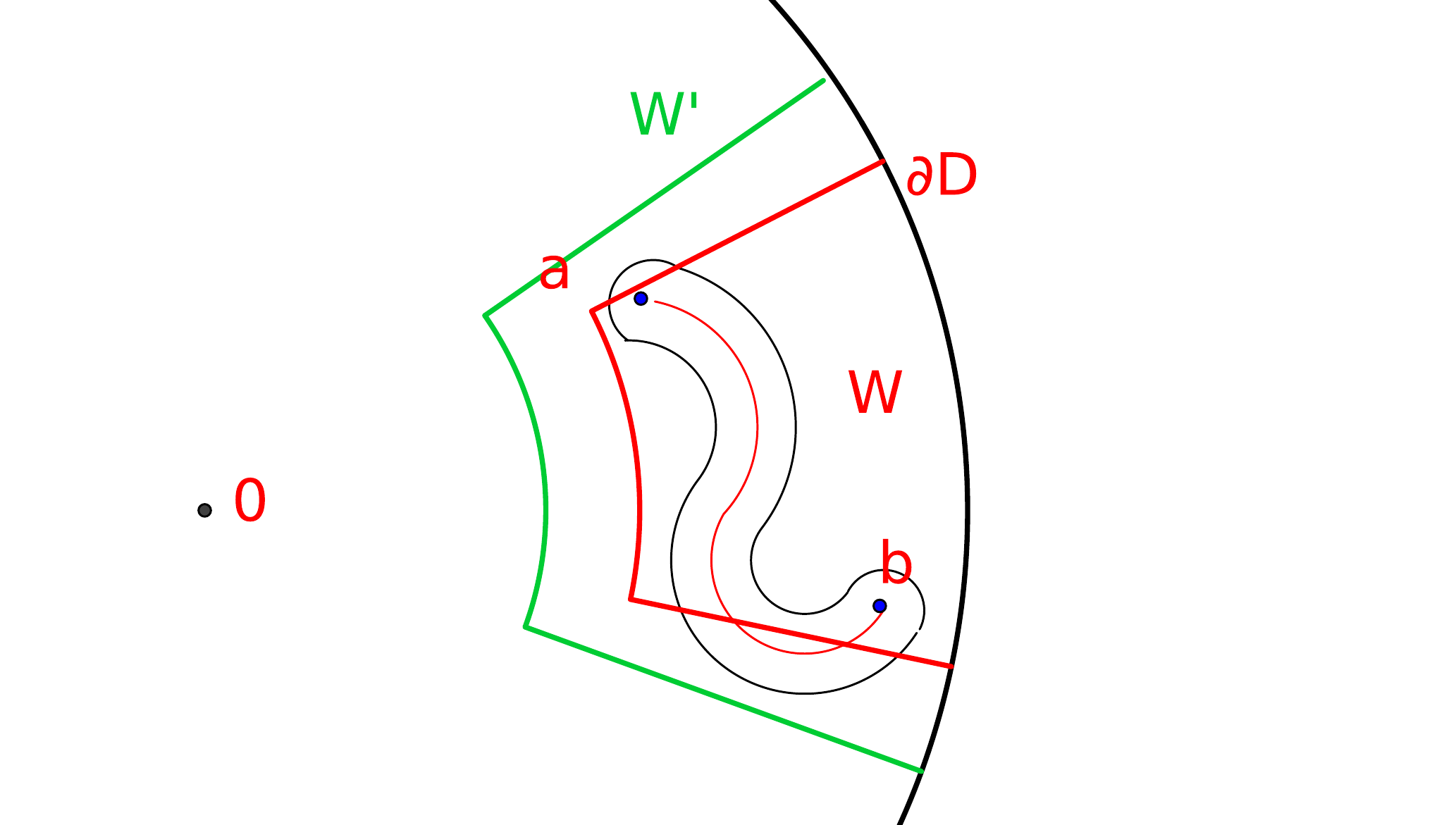}
\
 \par 
\quad \quad \quad 	The Carleson-H\"ormander lemma~\ref{ultraSepBall62} gives us
 : there is a  $\psi $  with  $0\leq \psi \leq 1$  so that  ${\rm{grad
 }}\psi $  is a Carleson measure in  ${\mathbb{D}}$  and\ \par 
\quad \quad \quad \quad \quad \quad  $\psi (z)=0$  when   $\ \left\vert{f(z)}\right\vert <\eta ^{\kappa
 },\ \psi (z)=1$  when  $\ \left\vert{f(z)}\right\vert >\eta .$ \ \par 
Let  $b:=\varphi (a),\ I_{a}:=\lbrack -s(1-\left\vert{a}\right\vert
 ),\ +s(1-\left\vert{a}\right\vert )\rbrack $  and parametrize
 the tube  $\displaystyle R(a,\varphi (a),s,\pi /2)$  :\ \par 
\quad \quad \quad \quad \quad \quad  $\displaystyle R(a,\varphi (a),s,\pi /2)=I_{a}{\times}\Gamma (a,b)$ \ \par 
then because  $\psi =0$  on  $I_{a}{\times}\lbrace a\rbrace $
  and  $\psi =1$  on  $I_{a}{\times}\lbrace b\rbrace ,$  by the
 known facts 1 and 2 and the construction of the  $R(a,\varphi
 (a),s,\pi /2),$  we have\ \par 
\quad \quad  $\displaystyle \forall t\in I_{a},\ 1=\psi (t,b)-\psi (t,a)=\int_{\Gamma
 (a,b)}{{\rm{grad }}\psi \cdot ds}\Rightarrow 1\leq \int_{\Gamma
 (a,b)}{\left\vert{{\rm{grad }}\psi }\right\vert ds}.$ \ \par 
Now we integrate with respect to  $t$  :\ \par 
\quad \quad \quad \quad \quad \quad  $\displaystyle 2s(1-\left\vert{a}\right\vert )=\int_{I_{a}}{1dt}\leq
 \int_{I_{a}}{\int_{a}^{b}{\left\vert{{\rm{grad }}\psi }\right\vert
 ds}dt}\leq \frac{\pi }{2}\int_{R(a,b,s,\pi /2)}{\left\vert{{\rm{grad
 }}\psi }\right\vert dm},$ \ \par 
because the Lebesgue measure on  $\displaystyle R(a,b,s,\pi /2)$
  is  $\pi /2$  equivalent to the product measure.\ \par 
\quad \quad \quad 	This gives the estimate for the points  $a$  in  $W$  such that
  $R(a,b,s,\pi /2)$  is in  $W'=W(\zeta ,\pi h/2)$  because\ \par 
\quad \quad \quad \quad \quad \quad  $\displaystyle \ \sum_{a\in F_{W}}{(1-\left\vert{a}\right\vert
 )}\leq \frac{1}{2s}\sum_{a\in F_{W}}{\frac{\pi }{2}\int_{R(a,b,...)}{\left\vert{{\rm{grad
 }}\psi (z)}\right\vert dm(z)}}\leq $ \ \par 
\quad \quad \quad \quad \quad \quad \quad \quad \quad \quad \quad \quad \quad \quad  $\displaystyle \leq \frac{\pi }{4s}\int_{W'}{\left\vert{{\rm{grad
 }}\psi (z)}\right\vert dm(z)}\leq \frac{\pi ^{2}}{8s}Ch,$ \ \par 
the tubes being disjoint by lemma~\ref{SepDisc610} and the last
 inequality because  $\displaystyle \ \left\vert{{\rm{grad }}\psi
 (z)}\right\vert dm(z)$  is a Carleson measure.\ \par 
\quad \quad \quad 	Hence the sequence  $A$  is Carleson and separated so it is
  $H^{\infty }({\mathbb{D}})$  interpolating.\ \par 
For the sequence  $B$  we proceed analogously and we get that
  $B$  is still separated and Carleson, hence  $H^{\infty }({\mathbb{D}})$
  interpolating. Because the union  $S=A\cup B$  is separated,
 we get that  $S$  is still  $H^{\infty }({\mathbb{D}})$  interpolating
 and this finishes the proof of the direct part of the theorem.\ \par 
\ \par 
\quad \quad \quad 	For the converse part of the theorem let  $S$  be an interpolating
 sequence for  $H^{\infty }({\mathbb{D}}).$  Then  $S$  is  $\delta
 $ -separated, hence discrete, so take any restricted or Hoffman
 partition  $(A,B)$  of  $S.$  Because  $S$  is  $H^{\infty }({\mathbb{D}})$
  interpolating, there is a  $f\in H^{\infty }({\mathbb{D}})$
  such that  $f=0$  on  $A,\ f=1$  on  $B$  and  $\ {\left\Vert{f}\right\Vert}_{\infty
 }\leq C.$  This means that  $g:=f/C$  ultra-separates the sequence
  $S$  for any  $\kappa >1,$  hence the theorem.  $\hfill\blacksquare $ \ \par 
\ \par 
\quad \quad \quad 	Now the answer to the question by A. Hartman :\ \par 
\begin{Corollary}
Let  $S=A\cup B,$  where  $(A,B)$  is a restricted or a Hoffman
 partition of  $S.$  Suppose that the Blaschke product  $B_{A},$
  precisely zero on  $A,$  verifies that  $\inf _{b\in B}\left\vert{B_{A}(b)}\right\vert
 \geq \eta >0,$  then  $S$  is a  $H^{\infty }({\mathbb{D}})$
  interpolating sequence.\par 
\end{Corollary}
\quad \quad \quad 	Proof.\ \par 
We have that  $0<\eta ^{\kappa },$  where  $\kappa $  is the
 Carleson constant, so we can apply the theorem~\ref{ultraSepBall61}.
  $\hfill\blacksquare $ \ \par 
\ \par 

\bibliographystyle{/usr/share/texmf/bibtex/bst/base/plain}

\end{document}